\input cyracc.def
\font\tencyr=wncyr10 at 12pt
\def\cyr{\tencyr\cyracc}
\font\twelvemsb=msbm10 at 12pt
\newfam\msbfam
\textfont\msbfam=\twelvemsb
\def\Bbb#1{\fam\msbfam\relax#1}
\documentstyle[12pt]{article}
\advance\textheight by \topskip
\def\cfrac#1#2{{\displaystyle #1 \over \displaystyle #2}}
\newtheorem{theorem}{Theorem}
\begin{document}
\title{Compact Jacobi matrices: \\
from Stieltjes to Krein and M$(a,b)$}
\author{Walter Van Assche \thanks{Senior Research Associate of the
Belgian National Fund for Scientific Research} \\
 Department of Mathematics \\
 Katholieke Universiteit Leuven \\
 Celestijnenlaan 200 B \\
 B-3001 Heverlee (Leuven) \\
 Belgium}
\maketitle
\begin{abstract}
In a note at the end of his paper {\it Recherches sur les fractions
continues}, Stieltjes gave a necessary and sufficient condition
when
a continued fraction is represented by a meromorphic function.
This result is related to the study of compact Jacobi matrices. We
indicate how this notion was developped and used since Stieltjes,
with special attention to the results by M.G. Krein.
We also pay attention to the perturbation of a constant Jacobi
matrix by a compact Jacobi matrix, work which basically started
with
Blumenthal in 1889 and which now is known as the theory for the
class
$M(a,b)$.
\end{abstract}

\section{A theorem by Stieltjes}

Stieltjes' research in {\it Recherches sur les fractions continues}
\cite{stieltjes} deals with continued fractions of the form
$$   \cfrac{1}{\alpha_1z + \cfrac{1}{\alpha_2 + \cfrac{1}{\alpha_3
z +
  \cdots + \cfrac{1}{\alpha_{2n} + \cfrac{1}{\alpha_{2n+1} z +
\ddots}}}}}\ , $$
where the coefficients $\alpha_k$ are real and positive. Such a
continued
fraction is nowadays known as an S-fraction, where the S stands for
Stieltjes.
By setting $b_0 = 1/\alpha_1$ and $b_n = 1/(\alpha_n \alpha_{n+1})$
for
$n \geq 1$, and by the change of variable $z=1/t$, this continued
fraction can be written as
\begin{equation} \label{eq:bfrac}
   \cfrac{b_0}{z + \cfrac{b_1}{1 + \cfrac{b_2}{ z +
    \cfrac{b_3}{1 + \cfrac{b_4}{ z + \ddots}}}}} =
    \cfrac{b_0}{1 + \cfrac{b_1t}{1 + \cfrac{b_2t}{ 1 +
    \cfrac{b_3t}{1 + \cfrac{b_4t}{ 1 + \ddots}}}}}\ ,
\end{equation}
where $b_k > 0$, which results from the positivity of the
$\alpha_k$.
Finally we can `contract' this fraction by using repeatedly
the identity
$$   z + \frac{a}{1+ b/c} = z + a - \frac{ab}{b+c} , $$
and the original S-fraction then changes to
\begin{equation} \label{eq:Jfrac}
  \cfrac{\lambda_0}{z+a_1 -\cfrac{\lambda_1}{z+a_2 -
\cfrac{\lambda_2}{
   z+ a_3 - \cfrac{\lambda_3}{z+a_4 - \ddots}}}} ,
\end{equation}
with
\begin{equation}  \label{eq:btoJ}
   a_n = b_{2n-2} + b_{2n-1}, \quad \lambda_n = b_{2n} b_{2n-1} .
\end{equation}
Such a continued fraction is known as a J-fraction, where the
letter J
stands for Jacobi. This J-fraction and the
original S-fraction are `nearly' equivalent in the sense that the
$n$th
convergent
of the J-fraction is identical to the $2n$th convergent of the
S-fraction.

During his work in \cite{stieltjes}, in particular the sections \S
68--69,
Stieltjes shows that the convergents of (\ref{eq:bfrac}) are given
by
$$  \frac{P_{n}(z)}{Q_n(z)} = b_0  t \ \frac{U_n(t)}{V_n(t)}, $$
where $U_n$ and $V_n$ are polynomials,
and the convergence of the series $\sum_{k=1}^\infty b_k$ is
necessary and
sufficient for the convergence
\begin{equation}  \label{eq:sepconv}
  \lim_{n \to \infty} U_n(t) = u(t), \quad
    \lim_{n \to \infty} V_n(t) = v(t),
\end{equation}
for every $t \in {\Bbb C}$, uniformly on compact sets. The
functions $u$ and
$v$ are thus both
entire functions as they are uniform limits of polynomials. Hence
the continued
fraction (\ref{eq:bfrac})
converges to $$  \lim_{n \to \infty} \frac{P_n(z)}{Q_n(z)} =
\frac{1}{\alpha_1z}
    \frac{u(1/z)}{v(1/z)} = F(z)  $$
and the function $F$ is meromorphic in the complex $t$-plane and
meromorphic
in the complex $z$-plane without the origin. Furthermore the zeros
of $U_n$ and
$V_n$ are all real (and negative) and they interlace (nowadays a
well known
property for orthogonal polynomials, observed a century ago by
Stieltjes), hence
$F$ has infinitely many poles in the $z$-plane, which accumulate at
zero.
Stieltjes then writes this function as
\begin{equation} \label{eq:Sint}
  F(z) = \frac{s_0}{\alpha_1 z} + \frac{1}{\alpha_1}
    \sum_{k=1}^\infty \frac{s_k}{z+r_k},
\end{equation}
where $\sum_{k=0}^\infty s_k = 1$ and $s_k > 0$ for every $k > 0$
$(s_0 \geq
0)$, and then uses the Stieltjes integral (which he introduced
precisely for
such purposes) to write it as
$$  F(z) = \int_0^\infty \frac{d\Phi(u)}{z+u}, $$
where $\Phi$ is a (discrete) distribution function with jumps of
size
$s_k/\alpha_1$
at the points $r_k$ ($k >0$), and also at the origin if $s_0 > 0$.
So Stieltjes has proved the following result in \cite[\S
68--69]{stieltjes}:

\begin{theorem}    \label{thm:tracecl}
Suppose that $b_k > 0$ for $k \in {\Bbb N}$. Then
\begin{equation}  \label{eq:tracecl}
    \sum_{k=1}^\infty b_k < \infty
\end{equation}
is a necessary and sufficient condition in order that
the continued fraction (\ref{eq:bfrac}) converges to
$$  F(z) =  \int_0^\infty \frac{d\Phi(u)}{z+u}  =
\frac{1}{\alpha_1z}
    \frac{u(1/z)}{v(1/z)}, $$
where $u$ and $v$ are entire functions and $F$ is meromorphic for
$z \in
{\Bbb C} \setminus \{ 0\}$.
\end{theorem}

In fact the condition (\ref{eq:tracecl}) gives the {\it separate
convergence}
(\ref{eq:sepconv})
of the numerator and denominator of the convergents of the
continued fraction
and allows to write $F$ as the ratio of two entire functions in
${\Bbb C} \setminus \{ 0\}$.

In a note at the end of his paper \cite{stieltjes}, Stieltjes wants
to find all
the cases such that the continued fraction (\ref{eq:bfrac})
converges to
a function $F$ meromorphic for $t \in {\Bbb C}$ (or $z \in {\Bbb C}
\setminus \{
0 \}$), and not only those for which one has the separate
convergence
(\ref{eq:sepconv}) as in the
case
when (\ref{eq:tracecl}) holds. Stieltjes still assumes the $b_k$ to
be positive.
In the note he proves the following extension of Theorem
\ref{thm:tracecl}:

\begin{theorem}    \label{thm:comp}
Suppose that $b_k > 0$ for $k \in {\Bbb N}$. Then
\begin{equation}  \label{eq:comp}
    \lim_{n \to \infty} b_n = 0
\end{equation}
is a necessary and sufficient condition in order that
 the continued fraction (\ref{eq:bfrac}) converges to
$$  F(z) =  \int_0^\infty \frac{d\Phi(u)}{z+u}
    $$
where  $F$ is meromorphic for $z \in {\Bbb C} \setminus \{ 0\}$.
\end{theorem}

He proves the necessity of the condition in \S 3 of the note, and
the
sufficiency in \S 4. Obviously condition (\ref{eq:tracecl}) implies
(\ref{eq:comp}), but the latter condition is weaker.

The condition (\ref{eq:comp})
given in Theorem \ref{thm:comp} is also sufficient in the case
where the coefficients $b_n$ are allowed to be complex. This result
was
proved by Van Vleck (see Section 6). But for complex $b_n$
condition
(\ref{eq:comp}) is no longer necessary, as was shown by Wall
\cite{wall1}
\cite{wall2}.

\section{Compact Jacobi operators}
For the J-fraction (\ref{eq:Jfrac}) the condition (\ref{eq:comp})
is equivalent
to
\begin{equation}  \label{eq:Jcomp}
   \lim_{n \to \infty} a_n = 0, \qquad \lim_{n \to \infty}
\lambda_n = 0.
\end{equation}
Furthermore, the stronger condition (\ref{eq:tracecl}) is
equivalent to
\begin{equation} \label{eq:Jtracecl}
   \sum_{n=0}^\infty (|a_{n+1}| + \sqrt{\lambda_n}) < \infty .
\end{equation}
Hence, Stieltjes' results in terms of the J-fraction
(\ref{eq:Jfrac})
show that the J-fraction converges to a meromorphic function
$F$ in ${\Bbb C} \setminus \{0\}$
if and only if (\ref{eq:Jcomp}) holds, and this meromorphic
function
is given by the ratio $(\alpha_1z)^{-1}u(1/z)/v(1/z)$, with $u$ and
$v$ entire
functions, if and only if (\ref{eq:Jtracecl}) holds. The
convergence
holds uniformly on compact sets of the complex plane excluding the
poles
of $F$, which accumulate at the origin.
In Stieltjes' analysis he always worked with S-fractions for which
$b_n > 0$ for all $n$, which gives certain restrictions to the
coefficients $a_n$ and $\lambda_n$ of the J-fraction, but in fact
the
results also hold for general real $a_n$ and positive $\lambda_n$.

With the coefficients of the J-fraction (\ref{eq:Jfrac}) one can
construct
an infinite tridiagonal {\it Jacobi matrix}
           $$  J = \left( \begin{array}{ccccc}
  -a_1             & \sqrt{\lambda_1} &   0       & 0         &
\cdots \\
  \sqrt{\lambda_1} &  -a_2     & \sqrt{\lambda_2} & 0         &
\cdots \\
   0        & \sqrt{\lambda_2} &  -a_3     & \sqrt{\lambda_3} &
\cdots \\
   0        &   0       & \sqrt{\lambda_3} & \ddots    & \ddots \\
                 0        &   0       &    0      & \ddots    &
                            \end{array} \right) . $$
With this infinite matrix we associate an operator, which we also
call $J$,
acting on the Hilbert space $\ell_2$ of square summable sequences.
If the coefficients $a_n$ and $\lambda_n$ are bounded, then this
operator
is a self-adjoint and bounded operator, which we call the {\it
Jacobi operator}.
In order to find eigenvalues and eigenvectors, one needs to solve
systems of the form $J u = x u$, where $u \in \ell_2$ and $x$ is an
eigenvalue, which, when it exists, will be real due to the
self-adjointness.
This readily leads to a three-term recurrence relation
\begin{equation}  \label{eq:threetold}
    x u_n = \sqrt{\lambda_n} u_{n-1} - a_{n+1} u_n +
\sqrt{\lambda_{n+1}}
u_{n+1}, \qquad  n \geq 0,
\end{equation}
where $u_{-1} = 0$. The solution when $u_0=1$ is such that
$u_n=p_n(x)$ is a polynomial of degree $n$ in the variable
$x$ and this is precisely the denominator polynomial for the $n$th
convergent of the J-fraction. Another solution, with $u_0=0$ and
$u_1=1$
gives a polynomial $u_n=p_{n-1}^{(1)}(x)$ of degree $n-1$, and this
is the
numerator polynomial for the $n$th convergent of the J-fraction.
Applying the spectral theorem to the Jacobi operator $J$ shows that
there is a positive measure $\mu$ on the real line such that $J$ is
unitarily isomorphic to the multiplication operator $M$ acting on
$L_2(\mu)$
in such a way that the  unit vector $e_0 = (1,0,0,\ldots) \in
\ell_2$
(which is a cyclic vector) is mapped to the constant function
$x \mapsto 1$, and $J^n e_0$ is mapped to the monomial
$x \mapsto x^n$.
A simple verification, using the three-term
recurrence relation, shows that the unitary isomorphy also maps the
$n$th
unit vector $e_n = (\underbrace{0,0,\ldots}_{n {\rm\
zeros}},1,0,0,\ldots) \in
\ell_2$ to the polynomial $p_n$, and since $\langle e_n, e_m
\rangle =
\delta_{m,n}$ in the Hilbert space $\ell_2$, this implies
that $\langle p_n,p_m \rangle = \int p_n(x) p_m(x) \, d\mu(x) =
\delta_{m,n}$
in the Hilbert space $L_2(\mu)$, showing that we are dealing with
orthogonal polynomials. For more regarding this connection between
spectral theory and orthogonal polynomials, see e.g.,
\cite{dombrowski},
\cite{mate}, \cite[\$ XII.10, pp.~1275--1276]{dunford},
\cite[pp.~530--614]{stone}. Unfortunately, the spectral theorem
(and the Riesz
representation theorem) came decades after Stieltjes so that
Stieltjes
was not using the terminology of orthogonal polynomials, even
though
he clearly was aware of this peculiar orthogonality property of the
denominator polynomials, as can be seen from \S 11 in
\cite{stieltjes}.

The spectrum of the operator $J$ corresponds to the support of the
{\it spectral measure} $\mu$. This spectrum is real since $J$ is
self-adjoint.
The measure $\mu$ in general consists of an absolutely continuous
part,
a singular continuous part and an atomic (or discrete) part, and
the
supports of these three parts correspond to the absolutely
continuous
spectrum, the singular continuous spectrum and the point spectrum.
The
point spectrum is the closure of the set of eigenvalues of $J$, and
thus
all $x \in {\Bbb R}$ for which $\sum_{n=0}^\infty p_n^2(x) <
\infty$
are in the point spectrum. Moreover, one can show that
$$   \mu(\{x\}) = \left( \sum_{n=0}^\infty p_n^2(x) \right)^{-1}.
$$

In terms of the Jacobi operator $J$, Stieltjes' results can be
formulated
as follows. The spectral measure $\mu$ corresponds to the
distribution
function $\Phi$ in (\ref{eq:Sint}) and
is purely atomic and the point
spectrum has $0$ as its only accumulation point if and only if
(\ref{eq:Jcomp}) holds. When (\ref{eq:Jtracecl}) holds, then
the eigenvalues are the reciprocals of the zeros of the entire
function $v$, which is obtained as the limit
$$   \lim_{n \to \infty} x^{-n} p_n(1/x) = v(x). $$
Since $v$ is the uniform limit (on compacta) of a sequence of
polynomials,
it follows that $v$ is an entire function. In fact, it is a
canonical
product completely determined by its zeros, and the order of this
canonical
product is less than or equal to one. Therefore
the sum $\sum_{n=0}^\infty 1/|z_k|$ converges, where $z_k$ are the
zeros
of the entire function $v$.

In modern terminology, the condition (\ref{eq:Jcomp}) implies that
the Jacobi
operator $J$ is a {\it compact operator}. In general, a linear
operator $A$
acting on a Hilbert space ${\cal H}$
is called compact if it maps the unit ball in ${\cal H}$ onto a set
whose
closure is compact.
In other words, $A$ is compact if for every bounded sequence
$\psi_n$
$(n \in {\Bbb N})$
of elements in the Hilbert space ${\cal H}$, there is always a
subsequence
 $A \psi_n$ $(n \in \Lambda \subset {\Bbb N})$ that converges.
Compact operators are sometimes also known as completely continuous
operators, but this terminology is not so much in use anymore.
It is not hard to see that for an operator associated with
a banded matrix $A$ of bandwidth $2m+1$, i.e., $A_{i,j} = 0$
whenever $|i-j|>m$,
for some fixed $m$, compactness is equivalent with the condition
that $\lim_{n \to \infty} A_{n,n+k} = 0$ for every $k$ with $-m
\leq k \leq m$
\cite[\S 31, pp.~93--94]{akhiezerglazman}. Indeed, a diagonal
matrix $(m=0)$
is compact if and only if the entries on the diagonal tend to 0.
A banded matrix is of the form $A = A_0 + \sum_{k=1}^m (V^*A_k +
B_k V)$,
where $V$ is the shift operator and $A_0, A_k,B_k$
$(k=1,2,\ldots,m)$
are diagonal operators, and since $V$ is bounded and compact
operators
form a closed two-sided ideal in the set of bounded operators,
this shows that $A$ is compact if and only if each of diagonal
matrices
$A_0, A_k, B_k$ is compact.
Hence a Jacobi operator is compact if and
only if (\ref{eq:Jcomp}) holds. The simplest linear operators are,
of course,
operators acting on a finite dimensional Hilbert space, in which
case
we are dealing with matrices. Next in degree of difficulty are the
compact operators, which can be considered as limits of finite
dimensional matrices. Indeed, the structure of the spectrum of a
compact
operator is quite similar to the spectrum of a matrix
since it is a pure point spectrum with only one accumulation point
at the origin, a result known as the {\it Riesz-Schauder theorem}.
This is in
perfect agreement with Stieltjes'
result (Theorem \ref{thm:comp}) and the poles of the meromorphic
function $F$
in fact correspond to the point spectrum (the eigenvalues) of the
operator
$J$. So, Stieltjes' theorem is an anticipation of the
Riesz-Schauder
theorem (proved by Schauder in 1930) regarding the spectrum of a
compact
operator, but restricted to tridiagonal operators.
Similarly, Stieltjes' Theorem \ref{thm:tracecl} is related to a
subclass
of the compact operators, namely those compact operators for which
$\sum_{n=0}^\infty |x_k| < \infty$, where $x_k$ are the eigenvalues
of the operators. These operators are known as {\it trace class
operators}.
One can show that a banded operator $A$ is trace class if
$\sum_{n=0}^\infty \sum_{k=-m}^m |A_{n,n+k}| < \infty$, hence the
condition
(\ref{eq:Jtracecl}) means that $J$ is trace class, in which case
the
eigenvalues are in $\ell_1$.

\section{Some orthogonal polynomials with compact Jacobi matrix}
Stieltjes' theorems were rediscovered half a century later during
the
investigation of (modified) Lommel polynomials. First, H. M.
Schwartz
\cite{schwartz}
considered continued fractions of the form (\ref{eq:bfrac}) but
allowed
the $b_k$ to be complex, and the more general J-fraction
(\ref{eq:Jfrac})
with complex $\lambda_k$ and $a_k$.

Later Dickinson \cite{dickinson}, Dickinson, Pollak, and Wannier
\cite{dickinsonpw}, and Goldberg \cite{goldberg} also considered
the polynomials $h_{n,\nu}$ satisfying the recurrence relation
$$   h_{n+1,\nu}(x) = 2x(n+\nu) h_{n,\nu}(x) - h_{n-1,\nu}(x), $$
with initial conditions $h_{-1,\nu}=0$ and $h_{0,\nu}=1$. These
polynomials
appear in the study of Bessel functions and allow to express
a Bessel function $J_{n+\nu}$ as a linear combination of two
Bessel functions $J_\nu$ and $J_{\nu-1}$ as
$$ J_{\nu+n}(x) = h_{n,\nu}(1/x) J_\nu(x) - h_{n-1,\nu+1}(1/x)
J_{\nu-1}(x), $$
reducing the investigation of the asymptotic behaviour of Bessel
function with high index to the investigation of the polynomials
$h_{n,\nu}$, which are known as Lommel polynomials.
considering $p_n = \sqrt{(n+\nu)/\nu}\ h_n$, the three-term
recurrence
is of the form
$$   x p_n(x) = \frac{1}{2\sqrt{(n+\nu)(n+\nu+1)}}\ p_{n+1}(x) +
               \frac{1}{2\sqrt{(n+\nu)(n+\nu-1)}}\ p_{n-1}(x), $$
which corresponds to a J-fraction and Jacobi operator with
coefficients
$a_n=0$ and $\lambda_n = [4(n+\nu)(n+\nu-1)]^{-1}$.
Clearly $\lim_{n \to \infty} \lambda_n =0$ so that Stieltjes'
Theorem \ref{thm:comp} holds, and we can conclude that the Lommel
polynomials are orthogonal with respect to an atomic measure
with support a denumerable set with accumulation point at the
origin.
The spectrum of the Jacobi operator can be identified completely
by investigating the asymptotic behaviour of the Lommel
polynomials,
and it turns out that the spectrum consists of the closure of the
set
$\{ 1/j_{k,\nu-1}: k \in {\Bbb Z} \}$, where $j_{k,\nu-1}$ are the
zeros of the Bessel function $J_{\nu-1}$. These points indeed
accumulate at the origin, but the origin itself is not an
eigenvalue
of the operator $J$.
Note that Goldberg \cite{goldberg} observed that the analysis of
Dickinson, Pollak, and Wannier \cite{dickinson} \cite{dickinsonpw}
was incomplete since they did not give any information whether or
not
the accumulation point $0$ had positive mass.
The Jacobi operator in this case is not trace class,
since (\ref{eq:Jtracecl}) is not valid. This is compatible with
the asymptotic behaviour $j_{n,\nu} \sim \pi n $ for the zeros
of the Bessel function.

For the Bessel functions there are several $q$-extensions, with
corresponding Lommel polynomials. For the Jackson $q$-Bessel
function
the $q$-Lommel polynomials were introduced by Ismail \cite{ismail}
who showed that these polynomials are orthogonal on a denumerable
set similar as for the Lommel polynomials but involving the
zeros of the Jackson $q$-Bessel function. For the Hahn-Exton
$q$-Bessel function the $q$-analogue of the Lommel polynomials
turn out to be Laurent polynomials and in \cite{koelink}
it is shown that they obey orthogonality with respect to
a moment functional acting on Laurent polynomials.

Other families of orthogonal polynomials with a compact Jacobi
matrix
include the Tricomi-Carlitz polynomials, for which the asymptotic
behaviour was recently studied by Goh and Wimp \cite{gohwimp}.
These polynomials satisfy the three-term recurrence relation
$$   (n+1)f_{n+1}(x) - (n+\alpha) x f_n(x) + f_{n-1}(x) = 0, $$
with $f_0 = 1$ and $f_{-1}=0$. For the orthonormal polynomials
$[n! (\alpha+n)/\alpha]^{1/2} f_n$ this gives
$a_n=0$ and $\lambda_n = n/[(n+\alpha)(n+\alpha-1)$,
so that $\lambda_n \to 0$ but the Jacobi operator is not trace
class.
The spectral measure now is supported on the set $\{ \pm
1/\sqrt{k+\alpha}:
k=0,1,2,\ldots \}$, which is indeed a denumerable set with
an accumulation point at the origin, and the elements are not
summable.
The Tricomi-Carlitz polynomials are also known as the
Carlitz-Karlin-McGregor polynomials \cite{asis} because Karlin and
McGregor
showed that they
turn out to be the orthogonal polynomials  for the imbedded random
walk of a queueing process with infinitely many servers and
identical
service time rates.
There are a number of other examples of orthogonal polynomials
arising from birth-and-death processes
for which the Jacobi operator is compact. Van Doorn \cite{doorn}
showed that the orthogonal polynomials for a queueing
process studied by B. Natvig in 1975, where potential customers
are discouraged by queue length, are orthogonal on a denumerable
set
accumulating at a point. The birth-and-death process
governing this queueing process has birth rates $\lambda_n =
\frac{\lambda}{n+1}$
($n \geq 0)$, which expresses that the rate of new customers
decreases
as the number $n$ of customers in the queue increases,
and death rates $\mu_0=0$ and $\mu_n = \mu$, which expresses that
the service
time does not depend on the queue length.
The corresponding orthogonal  polynomials then  satisfy the
three-term
recurrence relation
$$  - x Q_n(x) = \lambda_n Q_{n+1}(x) - (\lambda_n+\mu_n) Q_n(x)
   + \mu_n Q_{n-1}(x). $$
The orthonormal polynomials $q_n$ then satisfy
$$   x q_n(x) = \sqrt{\lambda_{n} \mu_{n+1}} q_{n+1}(x) +
    (\lambda_n+\mu_n) q_n(x) + \sqrt{\lambda_{n-1}\mu_n}
q_{n-1}(x), $$
and since
$$  \lim_{n \to \infty} \lambda_{n-1} \mu_n = 0 , \qquad
    \lim_{n \to \infty} \lambda_n+\mu_n = \mu, $$
it follows that these  polynomials correspond to a Jacobi
matrix $J$ which can be written as $J=\mu I + J_p$, where $J_p$
is a compact operator. Hence the orthogonality measure is
denumerable
with only one accumulation point at $\mu$. Van Doorn gives a
complete
description of the support of the accumulation point.
Chihara and Ismail \cite{chihism} studied these
polynomials in more detail and showed that the point $\mu$ is not
a mass point of the orthogonality measure, even though it is an
accumulation point of mass points. Chihara and Ismail also study
the
queueing process with birth and death rates
$$  \lambda_n = \frac{\lambda}{n+a}, \quad \mu_{n+1} = \frac{\mu
(n+1)}{n+a},
\qquad n\geq 0, $$
for which the Jacobi operator is again of the form $J = \mu+J_p$
with $J_p$ a compact operator. The case $a=1$ corresponds
to the situation studied by Natvig and van Doorn.
Another way to model a queueing process where potential customers
are discouraged by queue length is to take
$$  \lambda_n = \nu q^n, \quad \mu_n = \mu(1-q^n), \qquad 0 < q <
1, $$
in which case the decrease is exponential. The corresponding
orthogonal polynomials turn out to be $q$-polynomials of Al-Salam
and Carlitz \cite[\S 10 on p.~195]{chihara}.

The orthogonal polynomials $U_n$ associated with the
Rogers-Ramanujan continued fraction \cite{alsalam}
$$  U_{n+1}(x) = x(1+aq^n) U_n(x)  - bq^{n-1} U_{n-1}(x), \qquad 0
< q < 1, $$
have a compact Jacobi operator, which in addition belongs to the
trace
class.
Several orthogonal polynomials of basic hypergeometric
type ($q$-polynomials) have a Jacobi matrix which is a compact
operator that
belongs to the trace class, so that Stieltjes' Theorem
\ref{thm:tracecl}
can be used to find the orthogonality relation for these
polynomials.
Often this orthogonality relation can be written in terms of the
$q$-integral
$$  \int_0^b f(t)\, d_qt = b(1-q) \sum_{n=0}^\infty f(bq^n) q^n, $$
and for $a < 0 < b$ this $q$-integral is given by
$$ \int_a^b f(t)\, d_qt = \int_0^b f(t)\, d_qt + \int_0^{-a}
f(-t)\, d_qt
 , $$
so that the support of the measure is the geometric lattice
$\{ aq^k,bq^k, k=0,1,2,\ldots\}$ which is denumerable and has $0$
as the
only accumulation point.
The orthogonal polynomials of this type are
the {\it big $q$-Jacobi polynomials}, the {\it big $q$-Laguerre
polynomials},
the {\it little $q$-Jacobi polynomials}, the {\it little
$q$-Laguerre
polynomials} (also known as the {\it Wall polynomials} \cite[\S 11
on
p.~198]{chihara}), the {\it alternative $q$-Charlier polynomials},
and the {\it Al-Salam--Carlitz polynomials}, which we already
mentioned
earlier. These polynomials, with references to the literature, can
be found in
\cite{koekoek}.

\section{Krein's theorem}
The most interesting extension of Stieltjes's Theorem
\ref{thm:comp} on compact
Jacobi operators was made by M. G. Krein \cite{krein}. He
considered operators
of the form $g(J)$, where $J$ is a Jacobi operator and $g$ a
polynomial.
It is not so hard to see that the matrix for the operator $g(J)$ is
banded and symmetric, and when $J$ is a bounded operator, then
$g(J)$ is also bounded. The bandwidth of $g(J)$ is $2m+1$ when
$g$ is a polynomial of degree $m$. In \cite{krein}, Krein first
shows that a banded operator $A$ with matrix $(a_{i,j})_{i,j\geq
0}$ is compact
if and only if $\lim_{i,j \to \infty} a_{i,j} = 0$. But his main
result
is
\begin{theorem}[Krein]
In order that the spectrum of $J$ consists of a bounded set
with accumulation points in $\{x_1,x_2,\ldots,x_m\}$,
it is necessary and sufficient that $J$ is a bounded operator and
$g(J)$ is a compact operator, where
$g(x)=(x-x_1)(x-x_2)\cdots(x-x_m)$.
\end{theorem}

The polynomial $g$ of lowest degree for which $g(J)$ is a compact
operator
is known as the {\it minimal polynomial}, and the zeros of the
minimal
polynomial correspond exactly to the accumulation points of the
spectrum
of $J$. Krein explicitly refers to Stieltjes' work, which is a
special case
where the minimal polynomial is the identity $g: x \mapsto x$ and
the
spectrum is a compact subset of $(-\infty,0]$ or $[0,\infty)$ if we
make a reflection through the origin. Krein
mentions a remark by N. I. Akhiezer  that, by changing
Stieltjes' reasoning somewhat, one may by his  method obtain the
result
for one accumulation point without the restriction that the
spectrum
is on the positive (or negative) real axis. However, Krein
finds it improbable that the result for $m>1$ accumulation points
could be proved by Stieltjes' method.

In terms of the corresponding orthogonal polynomials, Krein's
theorem
says that when $g(J)$ is a compact operator, then the polynomials
will be orthogonal with respect to a discrete measure $\mu$ and the
support
of this measure has accumulation points at the zeros of $g$. It is
not so difficult to prove that orthogonal polynomials can have at
most
one zero in an interval $[a,b]$ for which $\mu([a,b])=0$. This
means
that also the zeros of the orthogonal polynomials will cluster
around these zeros of $g$.

In terms of the continued fraction (\ref{eq:Jfrac}) Krein's result
implies that the continued fraction will converge to a function $F$
which is meromorphic in ${\Bbb C} \setminus \{x_1,x_2,\ldots,x_m\}$
and the poles of this meromorphic function accumulate at the zeros
of $g$.

Recently it has been shown \cite{duran} that Krein' theorem can be
restated
in terms of orthogonal matrix polynomials, where the polynomials
have matrix coefficients with matrices from ${\Bbb R}^{m \times
m}$.
Orthogonal matrix polynomials satisfy a three-term recurrence
relation with matrix coefficients, and with these matrix
recurrence coefficients one can form a block Jacobi matrix, which
defines a self-adjoint operator, but now one does not have a single
cyclic vector, but a set of $m$ cyclic vectors. Consequently, the
spectrum
in not simple and the spectral measure is a (positive definite)
$m\times
m$ matrix of measures $M = (\mu_{i,j})_{1\leq i,j\leq m}$. Starting
with an
ordinary Jacobi matrix,
the matrix $g(J)$ is banded and can be considered as a block Jacobi
matrix,
where the subdiagonals are triangular matrices. If $g(J)$ is
compact,
then by the Riesz-Schauder theorem the spectrum $\sigma(g(J))$
of $g(J)$ has only one accumulation
point at the origin, which means that the spectral matrix of
measures
is discrete and the support, which is the support of the
trace measure $\sum_{j=1}^m \mu_{j,j}$, has only one accumulation
point at the origin. The spectral matrix of measures for $g(J)$ is
connected
with the spectral measure for $\mu$ and in particular
$\sigma(J) \subset g^{-1}(\sigma(g(J)))$, and since $\sigma(g(J))$
has only one accumulation point at 0, it follows that the spectrum
$\sigma(J)$ of $J$
has accumulation points at $g^{-1}(0)$, which are the zeros of $g$.

\section{The class M$(a,b)$ and Blumenthal's theorem}
Compact Jacobi operators have also shown to be of great use in
studying
orthogonal polynomials on an interval.
In this section we will change notation and consider the three-term
recurrence relation
\begin{equation} \label{eq:threetnew}
  x p_n(x) = a_{n+1} p_{n+1}(x) + b_n p_n(x) + a_n p_{n-1}(x),
\end{equation}
so that $(-a_{n+1},\sqrt{\lambda_n})$ in (\ref{eq:threetold})
corresponds to
$(b_n,a_n)$ in (\ref{eq:threetnew}). This notation is more common
nowadays.
If we consider orthogonal polynomials
satisfying a three-term recurrence relation with constant
coefficients,
$$   x\tilde{p}_n(x) = \frac{a}2 \tilde{p}_{n+1}(x) + b
\tilde{p}_n(x) +
   \frac{a}2 \tilde{p}_{n-1}(x), $$
with initial values $\tilde{p}_0 = 1$ and $\tilde{p}_{-1} = 0$,
then these
polynomials are given by
$$   \tilde{p}_n(x) = U_n((x-b)/a), $$
where the $U_n$ are the {\it Chebyshev polynomials of the second
kind},
defined as
$$   U_n(x) = \frac{\sin(n+1) \theta}{\sin \theta}, \qquad x = \cos
\theta. $$
For these polynomials the orthogonality relation is
$$  \frac{2}{\pi} \int_{-1}^{1} U_n(x) U_m(x) \sqrt{1-x^2} \, dx =
   \delta_{m,n}, $$
which follows easily from the orthogonality of the trigonometric
system
$\{ \sin k \theta, k=1,2,3,\ldots\}$. Hence, by an affine
transformation,
the polynomials $\tilde{p}_n$ $(n \in {\Bbb N}\}$ obey the
orthogonality
conditions
$$  \frac{2}{\pi a^2}
   \int_{b-a}^{b+a} \tilde{p}_n(x)\tilde{p}_m(x)
\sqrt{a^2-(x-b)^2}\,
dx = \delta_{m,n}. $$
Hence these polynomials are orthogonal on the interval $[b-a,b+a]$
and
they will serve as a comparison system for a large class of
polynomials
for which the essential support is $[b-a,b+a]$.  A measure
$\mu$ on the real line can always be decomposed as
$\mu = \mu_{ac} + \mu_{sc} + \mu_{d}$, where $\mu_{ac}$ is
absolutely
continuous, $\mu_{sc}$ is singular and continuous, and
$\mu_{d}$ is discrete (or atomic). The essential support of $\mu$
corresponds to the support of $\mu_{ac}+\mu_{sc}$ together with the
accumulation points of the support of $\mu_d$. Hence, if a measure
has essential support equal to $[b-a,b+a]$ then $\mu$ can have
mass points outside $[b-2a,b+2a]$, but the accumulation points
should be on this
interval. The Jacobi operator for $\tilde{p}_n$ has a matrix
with constant values
$$   J_0 = \left( \begin{array}{ccccc}
             b  &  a  &  0  &  0  & \cdots \phantom{\ddots} \\
             a  &  b  &  a  &  0  & \cdots \phantom{\ddots} \\
             0  &  a  &  b  &  a  & \cdots \phantom{\ddots} \\
             0  &  0  &  a  & \ddots    & \ddots \\
             0  &  0  &  0  & \ddots    &
              \end{array} \right) . $$
If we perturb this operator by adding to it a compact Jacobi
operator $J_p$, so
that we obtain a Jacobi operator
$$  J = \left( \begin{array}{ccccc}
             b_0  &  a_1  &  0    &  0  & \cdots \phantom{\ddots}
\\
             a_1  &  b_1  &  a_2  &  0  & \cdots \phantom{\ddots}
\\
             0  &  a_2  &  b_2  &  a_3  & \cdots \phantom{\ddots}
\\
             0  &  0  &  a_3  & \ddots    & \ddots \\
             0  &  0  &  0  & \ddots    &
              \end{array} \right)  = J_0 + J_p ,  $$
then the Jacobi operator $J$ has entries for which
\begin{equation}  \label{eq:Mab}
    \lim_{n \to \infty} a_n = \frac{a}2, \qquad
    \lim_{n \to \infty} b_n = b ,
\end{equation}
and we say that $J$ is a compact perturbation of $J_0$. There is a
very
useful result regarding compact perturbations of operators, which
is quite
useful in the analysis of orthogonal polynomials \cite{mate}.

\begin{theorem}[H. Weyl]
Suppose $A$ is a bounded and self-adjoint operator and $C$ is a
compact
operator, then $A+C$ and $A$ have the same essential spectrum.
\end{theorem}

Applied to our analysis of orthogonal polynomials, this means that
the orthogonal polynomials corresponding with a Jacobi operator
$J=J_0+J_p$ , where $J_p$ is compact, have an essential spectrum on
$[b-a,b+a]$, hence the orthogonality measure $\mu$ for these
polynomials has support $[b-a,b+a] \cup E$, where $E$ is at most
denumerable with accumulation points only at $b\pm a$. Compact
perturbations of the operator $J_0$ occur quite often, and in 1979
Paul Nevai \cite{nevai} introduced the terminology M$(a,b)$ for the
class of
orthogonal polynomials for which (\ref{eq:Mab}) holds.
The investigation of the class M$(a,b)$, however, goes back almost
a century.
The first to consider this class was O. Blumenthal, a student of
Hilbert,
whose Inaugural Dissertation \cite{blumenthal} was devoted to this
class.
In his dissertation, Blumenthal proves the following result
regarding
the continued fraction (\ref{eq:bfrac}):

\begin{theorem}[Blumenthal]
Streben die Gr\"ossen $b_n$ den endlichen von $0$
ver\-schie\-den\-en limites:
\begin{equation}  \label{eq:blumenthal}
  \lim b_{2n} = \ell, \qquad \lim b_{2n+1} = \ell_1
\end{equation}
zu, so liegen innerhalb des ganzen Intervalles
$$  \{ -(2\sqrt{\ell\ell_1}+\ell+\ell_1) \leq z \leq
        2\sqrt{\ell\ell_1}-\ell-\ell_1  \}  $$
\"uberall dicht Nullstellen der Funktionen-Reihe $Q_{2n}$,
ausserhalb desselben n\"ahern sich die Nullstellen mit wachsendem
$n$
einder endlichen Zahl von Grenzpunkten.%
\footnote{If the $b_n$ converge to positive limites
$b_{2n} \to \ell,\ b_{2n+1} \to \ell_1$, then the zeros of the
sequence
of functions $Q_{2n}$ will be dense in the interval
interval $[-(2\sqrt{\ell\ell_1}+\ell+\ell_1),
2\sqrt{\ell\ell_1}-\ell-\ell_1]$,
outside of which the zeros for increasing $n$ will approach
a finite number of limit points.}
\end{theorem}
In terms of the J-fraction, the convergence in 
(\ref{eq:blumenthal})
is equivalent with
$$  \lim_{n\to\infty} b_n = -\ell-\ell_1=b, \qquad \lim_{n \to
\infty}
    \a_n = \sqrt{\ell\ell_1} =a/2, $$
which corresponds to the class M$(a,b)$, and Blumenthal's
conclusion
is that the zeros of the denominator polynomials (the orthogonal
polynomials)
are dense on the interval $[b-a,b+a]$ (the essential spectrum) and
that outside this interval the zeros converge to a finite number of
limit points. The latter statement, however, turns out not to be
correct, since
outside the interval $[b-a,b+a]$ there can be a denumerable number
of limit points of the zeros, which can only accumulate at the
endpoints
$b\pm a$, which means that outside $[b-a-\epsilon,b+a+\epsilon]$
there
are a finite number of limit points, and this is true for every
$\epsilon > 0$ (but not for $\epsilon=0$). Except for this,
Blumenthal's
theorem is really a beautiful result and a nice complement to
Stieltjes'
Theorem \ref{thm:comp} which deals with the special case
$\ell=\ell_1=0$.

Blumenthal's proof of the theorem was based on a result by
Poincar\'e
\cite{poincare} which describes the ratio asymptotic behaviour
of the solution of a finite order linear recurrence relation when
the
coefficients in the recurrence relation are convergent.

\begin{theorem}[Poincar\'e]
If in the recurrence relation
$$   y_{n+k} = \sum_{j=0}^{k-1} a_{j,n} y_{n+j} $$
the recurrence coefficients  have limits
$$  \lim_{n \to \infty} a_{j,n} = a_j, \qquad 0 \leq j < k, $$
and if the roots $\xi_i$ $(i=1,2,\ldots,k)$ of the characteristic
equation
$$    z^k = \sum_{j=0}^\infty a_j z^j   $$
all have different modulus, then either $y_n = 0$ for all $n \geq
n_0$
or there is a root $\xi_{\ell}$ of the characteristic equation
such that
$$  \lim_{n \to \infty} \frac{y_{n+1}}{y_n} = \xi_\ell. $$
\end{theorem}

For a nice and comprehensive proof, see \cite{matenev}.
The case relevant for the class M$(a,b)$ corresponds to the second
order recurrence relation (\ref{eq:threetnew})
for orthogonal polynomials corresponding to a Jacobi matrix $J$,
and
(\ref{eq:Mab}) expresses the fact that the recurrence coefficients
have limits. The characteristic equation then is
$$    2x z = a z^2 + 2b z + a $$
for which the roots are
$$  \xi_1 = \frac{x-b + \sqrt{(x-b)^2-a^2}}{a}, \quad
  \xi_2 = \frac{x-b - \sqrt{(x-b)^2-a^2}}{a}. $$
These two roots have equal modulus whenever
$(x-b)^2-a^2 \leq 0$, hence for $x \in [b-a,b+a]$, so that this
simple
observation already gives the important interval. Poincar\'e's
theorem
then shows that for $x \notin [b-a,b+a]$ the ratio
$p_{n+1}(x)/p_n(x)$ converges to one of the two roots of the
characteristic equation. Poincar\'e's theorem does not tell you
which root, but in the case of orthogonal polynomials, we know that
for
$x$ large enough the ratio $p_{n+1}(x)/p_n(x)$ behaves like $x$ as
$x \to \infty$, hence we need to choose the root with largest
modulus
whenever $x$ is large enough.  This asymptotic behaviour can then
be used to
obtain information about the set of limit points of the zeros of
the
orthogonal polynomials, which is how Blumenthal arrived at his
results.
For a contemporary approach, see \cite{mate}.

Blumenthal's result thus deals with compact perturbations on
Chebyshev polynomials. If more can be said of the (compact)
perturbation
operator $J-J_0$, then more can also be said of the spectral
measure
for the Jacobi operator $J$. If $J-J_0$ is a trace class operator,
i.e.,
$$  \sum_{k=1}^\infty \left( |a_k-\frac{a}2| + |b_k-b| \right)  <
\infty, $$
then there is a beautiful theorem by Kato and Rosenblum
\cite[Thm.~4.4 on
p.~540]{kato} that tells
something about the nature of the spectral measure on the essential
spectrum
\cite{dombrowski}.

\begin{theorem}[Kato-Rosenblum]
Suppose $A$ is a self-adjoint operator in a Hilbert space ${\cal
H}$ and
and $C$ is a trace
class operator in ${\cal H}$ and that $A+C$ is self-adjoint.
Then the absolutely continuous parts of $A$ and $A+C$ are unitarily
equivalent.
\end{theorem}

The spectral measure for the operator $J_0$ is absolutely
continuous on
$[b-a,b+a]$, hence the Kato-Rosenblum theorem implies that
the orthogonal polynomials corresponding to the Jacobi
operator $J$ are orthogonal with respect to a measure with an
absolutely
continuous part in $[b-a,b+a]$. The measure can still have
a discrete part outside $[b-a,b+a]$. For even more information
regarding
this absolutely continuous part, one needs an even stronger
condition such as \cite{wva}
$$  \sum_{k=1}^\infty k \left( |a_k-\frac{a}2| + |b_k-b| \right) 
< \infty, $$
in which case $\mu'(x) = g(x) (x-b -a)^{\pm 1/2} (x-b+a)^{\pm
1/2}$,
where $g$ is continuous and strictly positive on $[b-a,b+a]$.
Furthermore,
in this case the number of mass points outside $[b-a,b+a]$ is
finite
and the endpoints $b \pm a$ are not mass points.

\section{Van Vleck's results}
The class M$(a,b)$ received a lot of attention the past two
decades,
starting with Nevai in \cite{nevai} who introduced the terminology
and obtained various results. See \cite{wva} and the references
given there for a survey on the class M$(a,b)$. In the mean time
it has become clear that this class has already been studied in
detail
almost a century ago by Blumenthal (see previous section), but also
by Edward B. Van Vleck. He studied the class in terms of continued
fractions,
much in the spirit of Stieltjes who also studied the class
of compact operators in terms of continued fractions. In
\cite{vanvleck2}
\cite{vanvleck}
Van Vleck considers continued fractions as in (\ref{eq:bfrac})
for which the coefficients converge. He does not require the
restrictions $b_n > 0$ and allows the coefficients to be complex.

\begin{theorem}[Van Vleck]
If in the continued fraction
\begin{equation} \label{eq:vleckfrac}
 \cfrac{b_0}{1 + \cfrac{b_1t}{1 + \cfrac{b_2t}{ 1 +
    \cfrac{b_3t}{1 + \cfrac{b_4t}{ 1 + \ddots}}}}}
\end{equation}
one has
$\lim_{n \to \infty} b_n = b$, then the continued fraction will
converge
in ${\Bbb C}$ except
\begin{enumerate}
\item along the whole or part of a rectilinear cut from $-1/4b$ to
$\infty$
with an argument equal to that of the vector from the origin to
$-1/4b$,
\item possibly at certain isolated points $p_1,p_2,p_3,\ldots$
\end{enumerate}
The limit of the continued fraction is holomorphic in ${\Bbb C}
\setminus
[-1/4b,\infty)$ except at the points $p_1,p_2,\ldots$ which are
poles.
\end{theorem}

Van Vleck's proof is again based on Poincar\'e's theorem. Van Vleck
actually shows that the exceptional points can have accumulation
points on the cut. In case all the $b_n$ are positive these
exceptional
points can only accumulate at the point $-1/4b$. Van Vleck also
considers the corresponding J-fraction (\ref{eq:Jfrac}). This will
also be a continued fraction with converging coefficients, and if
$b_n \to b$, then obviously
$$  a_n = b_{2n-2}+b_{2n-1} \to 2b, \quad \lambda_n =
b_{2n}b_{2n-1} \to b^2. $$
The limit  of the continued fraction (\ref{eq:Jfrac}) is equal
to $F(1/z)$, where $F(t)$ is the limit of the continued fraction
(\ref{eq:vleckfrac}). Hence $F(1/z)$ is analytic in the complex
plane
cut along the segment $[-4b,0]$, except at the points $1/p_1,
1/p_2, \ldots$
which are poles. The cut $[-4b,0]$ is indeed the essential
spectrum, since this J-fraction is one that corresponds to the
class
M$(2b,-2b)$. Van Vleck also considers the limiting case $b \to 0$
and thus was able to generalize Stieltjes' Theorem \ref{thm:comp}
for complex coefficients, showing that (\ref{eq:comp}) is a
sufficient
condition
(but not necessary condition, see Wall \cite{wall1} \cite{wall2})
 for a continued fraction to converge to a meromorphic function.

\end{document}